\documentclass[12pt,reqno]{amsart}

\usepackage{epsfig}
\usepackage{varioref}
\usepackage{hyperref}
\usepackage{ccaption}
\usepackage{amsmath,amsfonts,amsthm,bm,bbold}
\usepackage{abraces}

\usepackage{amssymb}

\def\finedim{{\hfill\hbox{\enspace${ \square}$}} \smallskip}    
\def\sqr#1#2{{\vcenter{\vbox{\hrule height .#2pt
     \hbox{\vrule width .#2pt height#1pt \kern#1pt \vrule
     width .#2pt} \hrule height .#2pt}}}}
\def\square{\mathchoice\sqr54\sqr54\sqr{4.1}3\sqr{3.5}3}
\def\Dim{\acapo \hbox{{\sl Proof:}$\;\; $}}
\def\DDD{\acapo \hbox{{\sl Proof [of Theorem \ref{teo3}]:}$\;\; $}}
\def\acapo{$\null$ \par {\vskip  0mm plus .2mm } \noindent }

\newcommand{\Tr}{{\mathrm{Tr}}}

\newcommand{\R}{{\mathbb{R}}}

\newcommand{\N}{{\mathbb{N}}}

\newcommand{\bP}{{\mathbb{P}}}
\newcommand{\bff}{{\bf f}}

\newcommand{\bE}{{\mathbb{E}}}
\newcommand{\Ba}{{\mathcal{B}}}

\newcommand{\Hi}{{\mathcal{H}}}

\newcommand{\Fo}{{\mathcal{F}}}

\newcommand{\cC}{{\mathcal{C}}}

\newtheorem{Theorem}{Theorem}

\newtheorem{Remark}{Remark}

\newtheorem{definition}{Definition}

\newtheorem{proposition}{Proposition}
\newtheorem{corollary}{Corollary}

\newtheorem{lemma}{Lemma}

\begin{document}

\title[Notes on the Ogawa integrability]{Notes on the Ogawa integrability and a condition for convergence in the multidimensional case}

\author{N. Cangiotti and S. Mazzucchi}
\address{Department of Mathematics, University of Trento, 38123 Povo (TN), Italy}

\maketitle

\begin{abstract}
The Ogawa stochastic  integral is shortly reviewed and formulated in the framework of abstract Wiener spaces. The condition of universal Ogawa integrability   in the multidimensional case is investigated, proving that it cannot hold in general without the introduction of a ``renormalization term''. Explicit examples are provided. \\

\noindent {\it Key words: Ogawa integral, noncausal calculus, stochastic integration.} 
\bigskip

\noindent {\it AMS classification }: 60H05, 46G12. 

\end{abstract}

\vskip 1\baselineskip

\section{Introduction}
\label{Intro}

After the introduction  of stochastic integral in the 1940s due to Kiyoshi It\^o and the developments of the It\^o calculus in the succeeding years, particular interest has been devoted to the hypothesis of causality, which are fundamental in stochastic integration theory. In fact, the It\^o calculus relies upon concepts as  adapted processes, filtration, martingale, conditions that seems to be consistent with a sort of principle of causality in physics. Hence,  for many years, the stochastic problems arising in physical modelling (e.g. the phenomenon of diffusion) could be effectively formulated  using It\^o calculus. 
Furthermore, the theory of martingales underlying in the It\^o calculus provides a powerful tool. 

However, at the end of 1960s, the interest to construct a new stochastic theory independently  from causality conditions began to take hold. In this context, many approaches have been developed. In particular Anatoliy Skorokhod defined, in 1970s, the so-called Skorokhod integral \cite{Sko75} and introduced the \emph{anticipative calculus}. A few years later, in 1979, Shigeyoshi Ogawa independently introduced the so-called Ogawa integral and the corresponding \emph{noncausal calculus} \cite{OG1979}. In these notes we have studied the latter with the aim to generalize the conditions for Ogawa integrability in the multidimensional case. 

In the following we shall adopt  Ogawa's  recent notations \cite{OG2007, OG2017}. Let us  set a probability space $(\Omega,\mathcal{F},\mathbb{P})$ and let $(W_t)_{t\in [0,1]}$ be the standard Wiener process  with natural filtration $\{\mathcal{F}_t\}$.  We define $\bold{H}$ as the set of real valued functions $f:[0,1]\times \Omega\to \R$ which are  measurable  with respect $B_{[0,1]}\times \mathcal{F}$  and such that  the following condition holds:
\[
\mathbb{P}\left( \int_0^1 |f(t,\omega)|^2dt < \infty \right )=1.
\]



Given an orthonormal basis $\{\phi_n\}$ of the Hilbert space $L^2([0,1],dt)$, let us consider the following formal random series
\begin{equation}
\label{formalseries}
S_\phi(f)\equiv \sum_{n=1}^{\infty}(f,\phi_n)(\phi_n,\dot W)
\end{equation}
where $(f,\phi_n)=\int_0^1 f(t)\bar \phi_n(t)dt$ denotes the inner product in $L^2([0,1],dt)$ and $(\phi_n,\dot W):=
 \int_0^1 \phi_n(t)dW_t$.
Now we can define a noncausal stochastic integral, i.e. the Ogawa integral.

\begin{definition}
\label{phi-int}
A function $f\in \bold H$ is said to be $\phi$-integrable (i.e. integrable with respect to the basis $\{\phi_n\}$) if the random series \eqref{formalseries} converges in probability. In this case this sum  is denoted
$
\int_0^1 fd_{\phi}W_t
$
and it is called the {\em Ogawa integral of $f$ with respect the basis $\{\phi_n\}$}. A function integrable with respect the basis $\{\phi_n\}$ is called $\phi$-integrable. 
\end{definition}

In Def. \ref{phi-int} the orthonormal basis $\{\phi_n\}$ plays an important role. The requirement of the independence of the existence as well as the value of the sum \eqref{formalseries} from the basis $\{\phi_n\}$ leads naturally to the definition of  \emph{universal integrability}.

\begin{definition}\label{u-int}
Let $f \in \bold H$. If $f$ is integrable in the sense of Def. \ref{phi-int} with respect any orthonormal basis and the value of the integral does not depend on the basis, then the function is called \emph{universally integrable} (\emph{u-integrable}).
\end{definition}


A different way to characterize the Ogawa integral, which comes directly from the It\^o-Nisio theorem \cite{ItoNisio}, is the following. We can consider the sequence of approximated processes as follows
\[
W_n^{\phi}(t)=\sum_{i=1}^n \int_0^t \phi_i(s)ds   \int_0^1 \phi_i(s)dW_s.
\]
According to the It\^o-Nisio theorem we have that the sequence $\{ W_n^{\phi}\}$ converges uniformly in $t\in [0,1]$ to $W_t$ with probability 1. Hence, the  Ogawa integral can also be defined as the limit of a sequence of Stieltjes integrals. In fact the following holds.
\begin{proposition}
Let $f \in \bold{H}$; then $f$ is $\phi$-integrable if and only if the sequence
\[
\int_0^1 f dW_n^{\phi}(t)
\]
of Stieltjes integral converges in probability. In particular we get
\[
\lim_{n \to \infty} \int_0^1 f dW_n^{\phi}(t)= \int_0^1 f d_{\phi}W_t.
\]
\end{proposition}

It is important to introduce the definition of regularity of an orthonormal basis. 
\begin{definition}
\label{RegularBasis}
 An orthonormal basis $\{\phi_n\}$ in $L^2([0,1],dt)$ is called \emph{regular} if 
\[
\sup_n \|u_n\|_{L^2} < \infty,
\]
where
\[
u_n(t)=\sum_{i\le n}\phi_i(t)\int_0^t\phi_i(s) ds.
\]
\end{definition}

\begin{Remark}
Two examples of regular basis are trigonometric functions and Haar functions. 
\end{Remark}

\begin{Remark}
The existence of a non-regular basis was proved by Pietro Majer and Maria Elvira Mancino in \cite{MajMan}.
\end{Remark}

%
%
%
%

\begin{Remark}
The results concerning the integrability with respect regular bases and with respect any orthonormal basis were studied by Ogawa \cite{OG1984} and then, in the context of Malliavin calculus, by David Nualart and Moshe Zakai \cite{NZ86}.
\end{Remark}

There are many approaches to the noncausal stochastic calculus (see e.g. \cite{Nualart88}). The Ogawa integral was extensively studied  also in relation   with the Skorohod integral \cite{NZ86} and the Stratonovich integral \cite{NZ89}.  Definitions \ref{phi-int} and \ref{u-int} have been extended to the case  of random fields \cite{NZ88, OG1991, Delgado98};  however a detailed study of the case where the integrand function is $d$-dimensional (with $d\geq 2$) is still lacking. In the present paper we are going to show that in the multidimensional case   the condition of  universal integrability cannot be fulfilled even in rather simple cases. 
 In Section 2 we use a result due to Roald Ramer \cite{RAMER} to show that in the particular case where $f(t, \omega)=\bm{\alpha}(W_t(\omega))$, where $\bm{\alpha}:\R^d\to \R^d$ is a smooth vector field, the convergence of \eqref{formalseries} to a limit which is independent of the basis $\{\phi_n\}$  requires the introduction of a ``renormalization term''. In Section 3 we provide some interesting examples.

\section{A renormalization term for multidimensional Ogawa integral on abstract Wiener spaces}

In the following we are going to present an equivalent definition of Ogawa integral with respect to Wiener process in the framework  of abstract Wiener spaces \cite{Gro1, Gro2,Kuo}.\\
 Let $(\Hi, \langle\,,\,\rangle)$ be the Hilbert space of absolutely continuous paths $\gamma:[0,1]\to \R^d$ such that $\gamma(0)=0$ and $\dot \gamma\in L^2([0,1],dt)$ ($\dot \gamma $ denoting the  weak derivative of $\gamma$), endowed with the inner product
 $$\langle \gamma, \eta \rangle =\int_0^1 \dot\gamma (s) \cdot\dot \eta (s)ds, \qquad \gamma,\eta\in \Hi.$$
 Let $\|\,\|$  denote the $\Hi$-norm, namely $\|\gamma\|^2=\int_0^1\dot \gamma (s)\cdot \dot \gamma(s)ds$, $ \gamma\in \Hi$.\\
Let $C=C([0,1];\R^d)$ be the Banach space of continuous paths $\omega:[0,1]\to \R^d$, endowed with the $\sup$-norm $|\;|$ and let $\bP$ be the Wiener measure on the Borel $\sigma$-algebra $\Ba(C)$ of $C$.
Since for $\gamma\in \Hi$ we have $|\gamma|\leq \|\gamma\|$, $\Hi$ is densely embedded in $C$. Denoted with $C^*$ the topological dual of $C$, we have the following chain of dense inclusions:
\begin{equation}
C^*\subset \Hi\subset C.
\end{equation}
In the following, with an abuse of notation we shall denote $\langle \eta , \omega \rangle$ the dual pairing between  two elements $\eta\in C^*$ and $\omega \in C$.\par
The finitely additive standard Gaussian measure $\mu$ defined as 
$$\mu(\cC_{P,D})=\int_{D}\frac{e^{-\frac{\|x\|^2}{2}}}{(2\pi )^{n/2}}dx,$$ on the cylinder sets $\cC_{P,D}\subset \Hi$ of the form $$\cC_{P,D}:=\{\gamma \in \Hi\colon P\gamma\in D\},$$
for some finite dimensional projection operator $P:\Hi\to \Hi$ (where $\dim(P\Hi)=n$) and some Borel set $D\subset \Hi$, 
does not extend to a $\sigma$-additive  measure on the generated $\sigma$-algebra.  Defined the cylinder sets in $C$ as 
$$\tilde \cC_{\eta_1,...,\eta_n;E}:=\{\omega \in C\colon (\langle \eta_1,\omega\rangle,\dots,\langle \eta_n,\omega\rangle)\in E\},$$
for some $n\in \N$, $\eta_1, \dots, \eta_n\in C^*$ and $E$ a Borel set of $\R^n$, we have that the intersection $\tilde \cC_{\eta_1,\dots,\eta_n;E}\cap\Hi$ is a cylinder set in $\Hi$. According to the fundamental results by Leonard Gross \cite{Gro1,Gro2}, the finite additive measure $\tilde \mu$ on the cylinder sets of $C$, defined as
$$\tilde \mu(\tilde \cC_{\eta_1,\dots,\eta_n;E}):=\mu(\tilde \cC_{\eta_1,\dots,\eta_n;E}\cap H)$$
extends to a $\sigma$ additive Borel measure on $C$ that coincides with the standard  Wiener measure $\bP$ in such a way that  for any $\gamma\in \Hi$ such that $\gamma$  is an element of $C^*$ the following holds
$$\int e^{i\langle\gamma, \omega \rangle}d\bP(\omega)=e^{-\frac{1}{2}\|\gamma\|^2}.$$
 This allows, in particular, to define, for any $\eta \in C^*$, a centered Gaussian random variable $n_\eta$ on $(C,\Ba(C),\bP)$ given by $n_\eta (\omega ):=\langle\gamma, \omega \rangle$. In particular, for $\eta ,\gamma\in C^*$, the following holds
\begin{equation}
\bE[n_\eta n_\gamma]=\int_0^1\dot \eta(s)\cdot \dot \gamma (s)ds=\langle \eta, \gamma \rangle,
\end{equation}
which shows that the map $n:C^*\to L^2(C,\bP)$ can be extended, by the density of $C^*$ in $\Hi$, to an unitary operator $n:\Hi\to L^2(C,\bP)$.

It is remarkable that, if $\gamma \in \Hi$, the Gaussian random variable $n_\gamma$ can be identified with the Paley-Wiener integral of $\dot \gamma\in L^2([0,1])$, i.e. $n_\gamma(\omega)=\int_0^1\dot\gamma (s)dW(s)$.

Given an orthogonal projector $P:\Hi\to\Hi$ with finite dimensional range, i.e. of the form $P(\gamma)=\sum_{i=1}^n\langle \gamma, e_i\rangle e_i$, with $\{e_1,\dots,e_n\}\subset \Hi$ orthonormal vectors generating $P(\Hi)$, it is possible to define a random variable $\tilde P:C\to \Hi$ as $\tilde P(\omega)=\sum_{i=1}^nn_{e_i}(\omega) e_i$. \\ More generally, a function $F:\Hi\to E$ on $\Hi$ with values in a Banach space $E$ is said to admit a {\em stochastic extension} $\tilde F:C\to E$ if for any sequence  $\{P_n\}$ of finite dimensional orthogonal projectors $P_n:\Hi\to\Hi$ converging strongly to the identity operator $I$,  the sequence of  random variables $\{F\circ \tilde P_n\}$ converges in probability to a random variable $\tilde F$ on $C$  (and the limit  does not depend on the sequence $\{P_n\}$). For further information and examples about abstract Wiener spaces and stochastic extensions see, e.g., \cite{Kuo, Amour2015}.

In this framework, the definition of Ogawa integral can be  reformulated.  Let us consider the $d-$dimensional canonical Wiener process, where $(\Omega, \Fo)=(C,\Ba(C))$ and $W_t(\omega )=\omega (t)$, $\omega \in C$.  Let ${\bf f}:[0,1]\times C\to \R^d$ be a  function in $\bold H$.
For any orthonormal basis $\{\phi_n\}$ of $L^2([0, 1];\R^d)$ we can construct a corresponding orthonormal basis $\{e_n\}$ of $\Hi$ as $e_n(s)=\int_0^s\phi_n(u)du$. In fact the map $U:L^2([0, 1]; \R^d)\to \Hi$ defined by 
\begin{equation}\label{mapU}U(\phi)(s)=\int_0^s\phi(u)du, \qquad \phi\in L^2([0, 1]; \R^d),\end{equation}is unitary with inverse given by $U^{-1}(\gamma)=\dot \gamma$, $\gamma \in \Hi$.  The finite dimensional approximations of the formal series \eqref{formalseries} can be equivalently written as
\begin{eqnarray}& &\sum_{i=1}^{n} \int_0^1 f(t,\omega)\phi_i(t)dt\int_0^1 \phi_i(t)dW_t\nonumber\\
& &=\sum_{i=1}^{n}n_{e_i}(\omega) \int_0^1 f(t,\omega)\dot e_i(t)dt\nonumber\\
& &= \int_0^1 f(t,\omega)\cdot  \dot\gamma_n (\omega )(t)dt\label{seq-O-2}
\end{eqnarray}
where 
\begin{equation}
\label{gamman}\gamma_n(\omega) :=\tilde P_n(\omega)=\sum_{i=1}^n e_in_{e_i}(\omega), \qquad \omega\in C.
\end{equation}
We  can say that $f$ is $\phi$-integrable if the sequence \eqref{seq-O-2} converges in probability. Analogously $f$ is defined to be universally Ogawa integrable if the limit does not depend on the sequence $\phi_n$ (or, equivalently,  on the sequence $\{e_n\}$). \\

In the following we shall show that  in the case $d\geq 2$ the condition of universal integrability is too strong and cannot be fulfilled even in the simplest cases.

Let us consider a $C^1$ vector field $\bm{\alpha}:\R^d\to \R^d$ and let $\bff:[0,1]\times C\to  \R^d$ defined as $\bff(t,\omega ):=\bm{\alpha}(\omega (t))$, $t\in [0,1]$.
Given an orthonormal basis $\{e_n\}$ of $\Hi$, let us consider the sequence $\{g_n\}$ of real random variables on $(C, \Ba(C),\bP)$ defined as
\begin{equation}\label{n1}g_n(\omega):=\int_0^1\bm{\alpha}(\omega (t))\cdot \dot \gamma_n(\omega (t))dt, \qquad \omega \in C.
\end{equation}
where $\gamma_n$ is defined in \eqref{gamman}.
Considered the function  $G: C\to\Hi$ defined as
\begin{equation}\label{functionG}
G(\omega)(t)=\int _0^t \bm{\alpha}(\omega (s))ds, \qquad \omega \in C, \ t\in [0,1],
\end{equation}
the functions $\{g_n\}$ can be represented by  the following inner product 
\begin{equation}\label{gn2}g_n(\omega)=\langle G(\omega), \tilde P_n(\omega)\rangle.\end{equation}

For $\omega\in C$, let $DG(\omega)$ denote the Frechet differential of $G$ evaluated in $\omega$, given by:
\begin{equation}\label{diffG}
DG(\omega)(\gamma)_j(t)=\int_0^t\nabla \alpha_j(\omega(s))\cdot \gamma(s)ds,
\end{equation}
where $\gamma\in \Hi$, and $j=1,\dots,d$.\\

We can now state the main result.
\begin{Theorem}\label{teo3}
For any orthonormal basis $\{e_n\}$ of $\Hi$, the sequence of renormalized finite dimensional approximations of the Ogawa integral, namely the sequence of real random variables $\{h_n\}$ on $(C, \Ba(C), \bP)$ defined as
\begin{eqnarray}
h_n(\omega)&=&g_n(\omega)-r_n(\omega)\nonumber\\
&=& \langle G(\omega), \tilde P_n(\omega)\rangle -\sum_{i=1}^n\langle e_i,DG(\omega)e_i\rangle,\label{hn}
\end{eqnarray}
converges in $L^2(C, \bP)$ and the limit is independent on the orthonormal basis $\{e_n\}$.
\end{Theorem}

The proof relies upon the following lemmas.
\begin{lemma} \label{lemmaRA1}Let $f:\R^n\to \R^n$ be a $C^1$ map such that $|f|$ and $\|JF\|_2$ belong to $L^2(\R^n,\mu)$, with $\|JF\|_2$ denoting the Hilbert-Schmidt norm of the Jacobian of $f$ and $\mu$ is the standard centered  Gaussian measure on $\R^n$. 
Then 
\begin{equation}
\int_{\R^n}\left(f(x)\cdot x-\Tr(Jf(x))\right)^2d\mu (x)\leq \int_{\R^n}\left(|f(x)|^2+\|Jf(x)\|_2^2\right)^2d\mu (x),
\end{equation}
where $\Tr(Jf(x))$ is the trace of the Jacobian of $f$.
\end{lemma}
For a detailed proof  of Lemma \ref{lemmaRA1} see \cite{RAMER}, where also  the following definition is introduced. 

\begin{definition} A function $G:C\to C$ with $G(C)\subset \Hi$  is said to be {\em $\Hi$-differentiable} if for any $\omega\in C$ the function $G_\omega:\Hi\to \Hi$ defined as $G_\omega(\gamma)=G(\omega +\gamma)$, $\gamma\in \Hi$, is Frechet differentiable at the origin in $\Hi$. Its Frechet derivative, namely the linear operator $DG_\omega(0)\in L(\Hi;\Hi)$,  will be denoted with the symbol $DG(\omega) $ and called the {\em $\Hi$-derivative of $G$ at $\omega$.}
\end{definition}
\begin{lemma} \label{lemmaRA2}Let $G:C\to C$, with $G(C)\subset \Hi$,  be a $\Hi$-differentiable map such that for any $\omega \in C$ the $\Hi$-derivative $DG(\omega)\in L(\Hi,\Hi)$ is an Hilbert-Schmidt operator. Let us assume furthermore that the maps $\|G\|:C\to \R$ and $\|DG\|_2:C\to \R$, where  $\|DG(\omega)\|_2$ denotes the Hilbert-Schmidt norm of $DG(\omega)$, belong to $L^2(\Omega, \bP)$. Let $\{e_i\}$ be an orthonormal basis of $\Hi$ and let $\{P_n\}$ and $\{\tilde P_n\}$ be the sequence of finite dimensional projectors on the span of $e_1,\dots,e_n$ and their stochastic extensions respectively. Then the sequence of random variables $\{h_n\}$ defined as 
$$h_n(\omega):= \langle G(\omega), \tilde P_n(\omega)\rangle -\Tr (P_nDG(\omega)), \qquad \omega\in C,$$
converges in $L^2(C, \bP)$ and the limit does not depend on the basis $\{e_i\}$.
\end{lemma}
The proof of Lemma \ref{lemmaRA2} is a direct consequence of Lemma 4.2 in \cite{RAMER}.

\DDD
It is straightforward to verify that the map $G:C\to C$  defined by \eqref{functionG} is $\Hi$-differentiable and its $\Hi$-derivative $DG$ is given by \eqref{diffG}.
Furthermore, for any $\omega\in C$, the operator $DG(\omega)$ is Hilbert-Schmidt. Indeed  $DG(\omega):\Hi\to \Hi$ is unitary equivalent to the linear operator $T:L^2([0, 1];\R^d)\to L^2([0, 1];\R^d)$ defined as 
\begin{equation}\label{operatorT}T=U^{-1}\circ DG(\omega)\circ U, \end{equation} where $U:L^2([0, 1]; \R^d)\to \Hi$ is the unitary operator defined in \eqref{mapU}. By direct computation it is simple to see that $T$ is explicitly given in terms of a kernel $K\in L^2([0,1]\times[0,1])$, i.e. for $\phi\in L^2([0, 1]; \R^d)$ and $t\in [0,1]$,
\begin{equation}
\label{Tformula}
(T\phi)_j(t)=\int_0^1K_j(t,t')\cdot \phi(t')dt', \quad j=1,\dots,d,
\end{equation}
where $K_j(t,t')=\nabla \alpha_j(\omega (t))\chi_{[0,t]}(t')$, $t,t'\in [0,1]$. By formula 4.32 in \cite{Moretti}, the Hilbert-Schmidt norm of $T$ is equal to:
\begin{eqnarray*}
\|T\|_2^2&=&\int_{[0,1]\times[0,1]}|K(t,t')|^2dtdt'=\sum_{j=1}^d\int_0^1\int_0^1|\nabla \alpha_j(\omega (t))|^2\chi_{[0,t]}(t')dt dt':\\
&=&\sum_{j=1}^d\int _0^1t|\nabla \alpha_j(\omega (t))|^2dt\leq \sum_{j=1}^d \int_0^1|\nabla \alpha_j(\omega (t))|^2dt<\infty,
\end{eqnarray*}
where the boundedness of the last expression follows by the continuity of the maps $t\mapsto \nabla \alpha_j(\omega (t))$. 
By the unitary equivalence of $T$ and $DG(\omega)$, we get
$$\|DG(\omega)\|_2^2=\sum_{j=1}^d\int _0^1t|\nabla \alpha_j(\omega (t))|^2dt<\infty.$$
Moreover, we have that 
\begin{eqnarray*}
\bE[\|G\|^2]&=& \int_0^1\bE[|\bm{\alpha}(\omega (t)|^2 ]dt =\\
&=&\int_0^1\int_{R^d}|\bm{\alpha}(x)|^2 \frac{e^{-\frac{|x|^2}{2t}}}{(2\pi t)^{d/2}}dxdt<+\infty\\
\bE[\|DG\|_2^2]&\leq & \sum_{j=1}^d\int_0^1\bE[|\nabla \alpha_j(\omega (t))|^2]dt= \\
&=&\sum_{j=1}^d\int_0^1\int_{\R^d}|\nabla \alpha_j(x)|^2\frac{e^{-\frac{|x|^2}{2t}}}{(2\pi t)^{d/2}}dxdt<\infty.
\end{eqnarray*}
By Lemma \ref{lemmaRA2} the sequence of random variables $\{h_n\}$ given by 
$$h_n(\omega)=\langle G(\omega), \tilde P_n(\omega)\rangle-\Tr(P_nDG(\omega))$$
converges in $L^2(C,\bP)$ an the limit does not depend on the orthonormal basis $\{e_i\}$. Furthermore, by direct computation,  the ``renormalization term'' $\Tr(P_nDG(\omega))$ is given by
$$\Tr(P_nDG(\omega))=\sum_{i=1}^n \langle e_i,DG(\omega)e_i\rangle=\sum_{i=1}^n\int_0^1 \dot e_i(t)\cdot (e_i(t)\cdot \nabla)\bm{\alpha}(\omega(t))dt.$$

\finedim 

\begin{corollary}
For any orthonormal basis $\{e_n\}$ of $\Hi$, the sequence $h_n$ defined in Theorem \ref{teo3} converges in probability and the limit is independent of the basis $\{e_n\}$.
\end{corollary}
\section{Examples}
According to Theorem \ref{teo3}, the condition of existence of the limit in probability of the sequence of random variables $\{g_n\}$ defined in \eqref{n1}, i.e. the Ogawa integrability of the function $f\in \bold H$, with   $f(t,\omega ):=\bm{\alpha}(\omega (t))$, $t\in [0,1]$, with respect to the orthonormal  basis $\{\phi_n\}$ of $L^2([0,1];\R^d)$ (with $\phi_n=\dot e_n$) is equivalent to the existence of the limit in probability of the ``renormalization term'' $r_n(\omega)=\Tr [P_nDG(\omega)]$. 
Analogously, the universal Ogawa integrability of $f$ is equivalent to the convergence  in probability of $r_n$ to a limit which does not depend on the basis $\{e_n\}$ of $\Hi$. In particular, if the linear operator $DG(\omega)\in L(\Hi,\Hi)$ is not trace class, then the convergent of  sequence $\Tr [P_nDG(\omega)]$ is not guaranteed and, in general,   its value depends on the orthonormal basis $\{e_n\}$. We are going to show that this problem occurs even in very simple cases.

Let $d=2$ and $\bm{\alpha}:\R^2\to\R^2$ is a  linear vector field of the form 
\begin{equation}\label{vectoralpha}\bm{\alpha}(x,y)=(h_1x+k_1y,h_2x+k_2y).\end{equation} In this case  the map $G:C\to \Hi$ is given by
$$G(\omega)(t)=\left(h_1\int_0^t\omega_1(s)ds+k_1\int_0^t\omega_2(s)ds,h_2\int_0^t\omega_1(s)ds+k_2\int_0^t\omega_1(s)ds\right),$$ whrere $  \omega=(\omega_1,\omega_2)\in C$.
The $\Hi$-derivative $DG(\omega)$ for any $\omega \in C$ is the linear operator $DG:\Hi\to \Hi$ simply given by 
$$DG(\gamma) (t)=\left(h_1\int_0^t\gamma_1(s)ds+k_1\int_0^t\gamma_2(s)ds,h_2\int_0^t\gamma_1(s)ds+k_2\int_0^t\gamma_2(s)ds\right),$$ with $  \gamma=(\gamma_1,\gamma_2)\in \Hi.$
We can compute explicitly the spectrum of the self-adjoint operator $|DG|=\sqrt{DG^*DG}$. 
Indeed, setting for notational simplicity $L\equiv DG^*DG$ we have, for $\eta, \gamma \in \Hi$:
\begin{eqnarray*}
\langle \eta, L\gamma\rangle&=&\langle DG\eta, DG\gamma\rangle\\
&=&\int_0^1(\eta_1(t), \eta_2(t)) A (\gamma_1(t),\gamma_2(t)) ^Tdt, 
\end{eqnarray*}
with $$ A=\left(\begin{array}{ll} h_1^2+h_2^2 & h_1k_1+h_2k_2\\ h_1k_1+h_2k_2 & k_1^2+k_2^2  \end{array}\right).$$
Hence, for $\gamma \in \Hi$ the vector $L(\gamma)\in \Hi$ is given by
$$L(\gamma )(t)^T=-\int_0^t\int_1^sA\gamma (r)^Tdrds.$$
$L$ is a compact operator and has a discrete spectrum. By introducing in $\R^2$  an orthonormal basis $\{u_1,u_2\} $ of eigenvectors of the  symmetric matrix $A$, with corresponding eigenvalues $a_1, a_2\in \R^+$, the eigenvectors $\{\gamma_n\}$ of $L$ can be represented as linear combination of $u_1$ and $u_2$, namely $\gamma_{n}=\eta_{n,1}u_1+\eta_{n,2}u_2$, with $\eta_{n,j}:[0,1]\to \R$.
 The  components $\{\eta_{n,j}\}$ of the eigenvectors (with eigenvalues $\lambda$) are solutions of 
$$\left\{\begin{array}{l}
\lambda_{n,j}\ddot \eta_{n,j}+a_j\eta_{n,j}=0\\
\dot \eta_{n,j}(1)=0\\
\eta_{n,j}(0)=0
\end{array}
\right.\qquad j=1,2,$$
which yields in the non-trivial case where $a_j>0$ the solutions   $\lambda _{n,j}= \frac{4 a_j}{\pi^2(1+2n)^2}$, with corresponding eigenvectors $\gamma_{n,j}(t) =\sin\left(\left(\frac{\pi}{2}+n\pi\right)t\right)u_j $, where $j=1,2$. 
Hence, we can conclude that $\|DG\|=\sqrt L$ is not trace class and in general the limit of $r_n=\Tr(P_nDG)$ does not necessary exist and, if it exists, its value depends on the sequence of projectors $\{P_n\}$ or, equivalently, on the choice of the orthonormal basis $\{e_n\}$ of $\Hi$.
In the following we are going to investigate the value that the ``renormalization term''  assumes for different choices of the orthonormal basis $\{e_n\}$.\\

Let us consider $L^2([0, 1];\R^2)$ and the following orthonormal basis 
\begin{align*}
\{\psi_n\}:&=\left \{ (1,0) , (0,1), \sqrt{2}(\cos(2\pi n t),0),\sqrt{2}(\sin(2\pi n t),0), \right . \\
&\qquad \left .   \sqrt{2}(0,\cos(2\pi n t)),\sqrt{2}(0,\sin(2\pi n t)) \right \}=\\
&=\{\psi_{0,x}, \psi_{0,y},\psi_{n,1},\psi_{n,2},\psi_{n,3},\psi_{n,4} \},
\end{align*}
with $n \in \N \setminus \{0\}$.
Rewriting formula (\ref{Tformula}) explicitly, we can compute
\[
\langle \psi_n, T\psi_n \rangle =\int_0^1 \psi_n(t)\cdot \left ( \int_0^t \psi_n(s) ds \cdot \nabla \right ) \bm{\alpha} (\omega(t)) dt,
\]
where $\bm{\alpha}:\R^2 \to \R^2$ is given by \eqref{vectoralpha} and $T:L^2([0, 1];\R^2)\to L^2([0, 1];\R^2)$ is defined by \eqref{operatorT}.   
For the vectors of the form $\psi_{n,j}$ with $j=1, \dots,4$ we have:
\begin{equation*}
\langle \psi_{n,j}, T\psi_{n,j} \rangle =  0;
\end{equation*}
while  for the two constant vectors
\begin{align*}
\langle \psi_{0,x}, T\psi_{0,x} \rangle&=\frac{h_1}{2};\\
\langle \psi_{0,y}, T\psi_{0,y} \rangle&=\frac{k_2}{2}.
\end{align*}
These results lead us to the following ``renormalization term" depending on the divergence of $\alpha$ for such basis:
\[
r_n=\Tr(P_nDG)=\sum_{i=1}^n \langle \psi_i, T \psi_i \rangle=\frac{1}{2}\nabla \cdot \bm{\alpha}.
\]

Let us now consider a different basis in $L^2([0, 1];\R^2)$:
\begin{align*}
\{\xi_n\}:&=\left \{   (1,0) , (0,1),(\cos(2\pi n t),\sin(2\pi n t)),(\sin(2\pi n t),\cos(2\pi n t)), \right . \\
&\qquad \left .  (-\cos(2\pi n t),\sin(2\pi n t)),(-\sin(2\pi n t),\cos(2\pi n t)) \right \}=\\
&= \{\xi_{0,x}, \xi_{0,y},\xi_{n,1},\xi_{n,2},\xi_{n,3},\xi_{n,4}\},
\end{align*}
with $n \in \N \setminus \{0\}$.
We use the same argument as before for the vectors 
\[
\xi_{n,1}=(\cos(2\pi n t),\sin(2\pi n t)).
\]

We obtain:
\begin{align*}
\langle \xi_{n,1},T\xi_{n,1}\rangle &=\int_0^1  \left( k_1 \frac{\sin ^2(\pi  n t) \cos (2 \pi  n t) }{\pi  n}+ h_1 \frac{  \sin (2 \pi  n t) \cos (2 \pi  n t)}{2 \pi  n}+  \right .\\
&\qquad \left. +  k_2 \frac{\sin (2 \pi  n t) \sin ^2(\pi  n t)}{\pi  n}+
h_2\frac{\sin ^2(2 \pi  n t) }{2 \pi  n} \right) dt=\\&=\frac{h_2-k_1}{4n\pi}=\frac{\nabla \times \bm{\alpha}}{4n\pi}.
\end{align*}

Analogously
\begin{align*}
\langle \xi_{n,2},T\xi_{n,2}\rangle &=\int_0^1  \left( h_2 \frac{\sin ^2(\pi  n t) \cos (2 \pi  n t) }{\pi  n}+ k_2 \frac{  \sin (2 \pi  n t) \cos (2 \pi  n t)}{2 \pi  n}+  \right .\\
&\qquad \left. +  h_1 \frac{\sin (2 \pi  n t) \sin ^2(\pi  n t)}{\pi  n}+
k_1\frac{\sin ^2(2 \pi  n t) }{2 \pi  n} \right) dt=\\&=\frac{k_1-h_2}{4n\pi}=-\frac{\nabla \times \bm{\alpha}(\omega(t))}{4n\pi}.
\end{align*}
and 
$$\langle \xi_{n,3},T\xi_{n,3}\rangle=\frac{k_1-h_2}{4\pi n}=-\frac{\nabla \times \bm{\alpha}}{4\pi n}$$
$$\langle \xi_{n,4},T\xi_{n,4}\rangle=\frac{-k_1+h_2}{4\pi n}=\frac{\nabla \times \bm{\alpha}}{4\pi n}$$

In this case the series $\sum_{i=1}^n \langle \xi_i, T \xi_i \rangle$ cannot converge absolutely and and the value of the ``renormalization term''  depends on the order of the terms in the sum. \\
%

At last we consider in the Hilbert space $\Hi$ the sequence of orthogonal projection operators onto the finite dimensional subspaces $H_n$ of piecewise linear paths of the form
\begin{equation}\label{plp}\gamma (t)=\sum_{i=0}^{n-1}  \mathbb{1}_{\left[ \frac{i}{n},\frac{i+1}{n} \right]}(t)\left(\gamma (i/n)+n\left(\gamma (i+1/n)-\gamma (i/n)\right)(t-i/n)\right), 
\end{equation}
with $t\in [0,1]$. An orthonormal basis of $H_n$ is provided, e.g., by the vectors 
$$\{(z_{n,i},0),(0,z_{n,i})\}_{i=0,...,n-1},$$
where
\[
z_{n,i}(t)=\sqrt{n}\chi_{ \left[ \frac{i}{n},\frac{i+1}{n} \right]}(t) \left ( t-\frac{i}{n}\right )+\frac{1}{\sqrt{n}}\chi_{ \left[ \frac{i+1}{n},1 \right]}(t),
\]
with $i=0,\dots,n-1$. We also notice that:
\[
\dot z_{n,i}(t)=\sqrt{n} \mathbb{1}_{\left[ \frac{i}{n},\frac{i+1}{n} \right]}(t).
\]
It is not difficult to compute
\[
\langle( \zeta_{n,i},0),DG(\zeta_{n,i},0)\rangle = \frac{h_1}{2n}, \quad \langle( 0,\zeta_{n,i}),DG(0,\zeta_{n,i})\rangle = \frac{k_2}{2n}, 
\]
Thereby we get
\begin{equation}\label{traccia-p-l}
\lim_{n\to\infty}\Tr[P_nDG]=\frac{1}{2}\nabla \cdot \bm{\alpha} .
\end{equation}
Eventually, the case where $\{P_n\} $ are the projectors on the subspaces of piecewise linear path described above, the study of the limiting behavior of the sequences $\{g_n\}$ and $\{r_n\}$ (defined respectively by \eqref{gn2} and \eqref{hn}) allows to compute explicitly the limit of the sequence $\{h_n\}$ for linear vector fields $\mathbb{\alpha}$. By Theorem \eqref{teo3} this limit is independent on the sequence of projectors.
\begin{lemma}
\label{Lemma3}
Let $G: C\to \Hi $ be a linear operator such that its restriction $G_\Hi$ on $\Hi$ is Hilbert-Schmidt. Let $\{P_n\}$ be a sequence of finite dimensional projection operators in $\Hi$ converging strongly to the identity. Then the sequences of random variables $\{g_n\}$ and $\{g'_n\}$ defined as:
\begin{eqnarray*}
g_n(\omega)&=&\langle G(\omega), \tilde P_n(\omega)\rangle, \qquad \omega \in C\\
g'_n(\omega)&=&\langle G(\tilde P_n(\omega)), \tilde P_n(\omega)\rangle, \qquad \omega \in C
\end{eqnarray*}
satisfy
\begin{equation}\label{conv-eq}
\lim_{n\to \infty}\bE[|g_n-g'_n|^2]=0.
\end{equation}

\end{lemma}
\Dim
\begin{eqnarray*}
\bE[|g_n-g'_n|^2]&=&\int |\langle G(\omega)-G(\tilde P_n(\omega)), \tilde P_n(\omega)\rangle|^2d\bP(\omega)\\
&=&\int |\langle G(\omega-\tilde P_n(\omega)), \tilde P_n(\omega)\rangle|^2d\bP(\omega)\\
&=&\int |\langle G(\sum_{j=n+1}^\infty e_j n_{e_j}(\omega)),\sum_{i=1}^n e_i n_{e_i}(\omega)\rangle|^2d\bP(\omega)\\
&=&\sum_{j,j'=n+1}^\infty \sum_{i,i'=1}^n \langle Ge_j,e_i\rangle \langle Ge_{j'},e_{i'}\rangle \bE [n_{e_j}n_{e_{j'}}n_{e_i}n_{e_{i'}}]\\
&=&\sum_{j=n+1}^\infty \sum_{i=1}^n( \langle Ge_j,e_i\rangle)^2 \\
&=&\sum_{j=n+1}^\infty \langle P_nGe_j,P_nGe_j\rangle,
\end{eqnarray*}
where in the third step we have applied It\^o-Nisio theorem.
By using the assumption that $G_\Hi$ is an Hilbert-Schmidt operator we obtain \eqref{conv-eq}.

\finedim
\begin{Theorem} \label{teoStra}Let $\bm{\alpha}$ be the linear vector field given by \eqref{vectoralpha} and $G:C\to \Hi$ the linear operator  \eqref{functionG}. Then the sequence of random variables $\{g_n\}$ defined by
$$g_n(\omega)=\langle G(\omega), \tilde P_n(\omega)\rangle, \qquad \omega \in C ,$$
where $\{P_n\}$ is the sequence of orthogonal projectors onto the subspaces $H_n$ of piecewise linear paths \eqref{plp}, converges in $L^2(C,\bP)$ to the Stratonovich integral
$$ \int_0^1 \bm{\alpha}(\omega(t))\circ d\omega(t).$$
\end{Theorem}
\Dim 
By lemma \ref{Lemma3} the sequence $\{g_n\}$ has the same limit of the sequence $\{g'_n\}$, where
$$g'_n(\omega)=\langle G(\tilde P_n(\omega)), \tilde P_n(\omega)\rangle, \qquad \omega \in C$$
if such a limit exists. Moreover the random variables $\{g'_n\}$ assume the following form
$$ g'n(\omega)=\int_0^1\bm{\alpha}(\omega_n(t))\cdot \dot \omega_n(t)dt, $$
where $\omega_n=\tilde P_n\omega \in \Hi$. By Wong-Zakai approximations results \cite{IkeWa}, , in the case where $\{P_n\}$  are projectors on piecewise linear paths,  the sequence $\{g'_n$ converges in $L^2(C,\mathbb{P})$ to the Stratonovich integral $ \int_0^1 \bm{\alpha}(\omega(t))\circ d\omega(t).$
\finedim

\begin{Theorem}
Let $\bm{\alpha}$ be the linear vector field given by \eqref{vectoralpha} and $G:C\to \Hi$ the linear operator  \eqref{functionG}. Then the sequence of random variables $\{h_n\}$ defined in Theorem \ref{teo3}, namely
$$h_n(\omega)=g_n(\omega)-r_n,$$
with $r_n=\Tr[P_nDG]$, converges to the It\^o integral.
$$\int_0^1 \bm{\alpha}(\omega(t))d\omega(t)$$ and the limit does not depend on the sequence $\{P_n\}$
\end{Theorem}
\Dim
By Theorem \ref{teo3} the sequence $\{h_n\}$ converges in $L^2(C,\mathbb{P})$ and the limit is independent of $\{P_n\}$. \\
In the case where $\{P_n\}$ are projectors onto subspaces of piecewise linear paths, we can compute explicitly the limit of both $\{g_n\} $ and  $\{r_n\}$.
Indeed, by Theorem \ref{teoStra} and formula \eqref{traccia-p-l}, we obtain
\begin{eqnarray*}
\lim_{n\to \infty }h_n(\omega)&=&\lim_{n\to \infty }g_n(\omega)-\lim_{n\to \infty }r_n\\
&=&\int_0^1 \bm{\alpha}(\omega(t))\circ d\omega(t)-\frac{1}{2}\nabla \cdot \bm{\alpha}
\end{eqnarray*}
where the limits are meant in $L^2(C,\mathbb{P})$. By 
the conversion formula between It\^o and Stratonovich integral 
\begin{equation}
\label{ItoStra}
\int_0^1 \bm{\alpha}(\omega(t))\circ d\omega(t)=\int_0^1 \bm{\alpha}(\omega(t))d\omega(t)+\frac{1}{2}\int_0^1\nabla \cdot \bm{\alpha}(\omega(t))dt,
\end{equation}
we obtain the final result

\finedim

\section{Acknowledgements}
Fruitful discussions with S. Albeverio, S. Bonaccorsi, V. Moretti and L. Tubaro are gratefully acknowledged. The first author acknowledges support by Riemann International School of Mathematics.

\end{document}